\documentclass{article}
\usepackage{amsmath,amsfonts,graphicx,url}

\newcommand{\abs}[1]{\lvert #1 \rvert}
\newcommand{\abss}[1]{\left\lvert #1 \right\rvert}
\date{}

\begin{document}

\title{On Gauss's First Proof of\\the Fundamental Theorem of Algebra}
\markright{The Fundamental Theorem of Algebra}
\author{Soham Basu and Daniel J.\ Velleman}

\maketitle

\begin{abstract}
Carl Friedrich Gauss is often given credit for providing the first correct proof of the fundamental theorem of algebra in his 1799 doctoral dissertation.  However, Gauss's proof contained a significant gap.  In this paper, we give an elementary way of filling the gap in Gauss's proof.
\end{abstract}

\section{Introduction.}

The fundamental theorem of algebra is the statement that every nonconstant polynomial with complex coefficients has a root in the complex plane.  According to John Stillwell \cite[pp.\ 285--286]{stillwell}:

\begin{quote}
It has often been said that attempts to prove the fundamental theorem began with d'Alembert (1746), and that the first satisfactory proof was given by Gauss (1799).  This opinion should not be accepted without question, as the source of it is Gauss himself.  Gauss (1799) gave a critique of proofs from d'Alembert on, showing that they all had serious weaknesses, then offered a proof of his own.  His intention was to convince readers that the new proof was the first valid one, even though it used one unproved assumption \ldots.  The opinion as to which of two incomplete proofs is more convincing can of course change with time, and I believe that Gauss (1799) might be judged differently today.  We can now fill the gaps in d'Alembert (1746) by appeal to standard methods and theorems, whereas there is still no easy way to fill the gap in Gauss (1799).
\end{quote}

\noindent Our goal in this paper is to respond to the challenge in Stillwell's final sentence by providing an elementary way to fill the gap in Gauss's 1799 proof \cite{gauss} of the fundamental theorem of algebra.

\section{Gauss's proof.}

In his 1799 proof, written when he was 22,
%the first of four proofs he gave, 
Gauss proved the fundamental theorem only for polynomials with real coefficients.  It is well known that this suffices to establish the theorem for all polynomials with complex coefficients.  To see why this is true, suppose the theorem holds for polynomials with real coefficients, and let $f(z) = c_N z^N + c_{N-1} z^{N-1} + \cdots + c_0$ be a nonconstant polynomial with complex coefficients.  Let $\overline{f}(z) = \overline{c_N} z^N + \overline{c_{N-1}} z^{N-1} + \cdots + \overline{c_0}$ be the polynomial whose coefficients are the complex conjugates of the coefficients of $f$, and let $g(z) = f(z)\overline{f}(z) = f(z) \overline{f(\overline{z})}$.  Then $g$ is a nonconstant polynomial with real coefficients, so by assumption it has a root $z_0$.  This means that $g(z_0) = f(z_0) \overline{f(\overline{z_0})} = 0$, so either $z_0$ or $\overline{z_0}$ is a root of $f$.

Given a polynomial $f$ of degree $N>0$ with real coefficients, Gauss considered the algebraic curves in the plane defined by the equations $\text{Re}(f(z)) = 0$ and $\text{Im}(f(z)) = 0$.  Each of these curves  consists of several continuous branches.  He showed that for sufficiently large $r$, each curve intersects the circle $\abs{z} = r$ at $2N$ points, and these intersection points are interleaved: between any two intersection points for one curve there is an intersection point for the other.  Gauss then claimed, without proof, that if a branch of an algebraic curve enters the disk $\abs{z} \le r$, then it must leave again.  Applying this fact to the curves $\text{Re}(f(z)) = 0$ and $\text{Im}(f(z)) = 0$, he concluded that if we start at one of the $2N$ intersection points of one of these curves with the boundary of the disk and follow the corresponding branch of the curve into the interior of the disk, then it must eventually emerge at one of the other intersection points.  Using this fact, together with the way the intersection points are interleaved, Gauss then argued that the two curves $\text{Re}(f(z)) = 0$ and $\text{Im}(f(z)) = 0$ must intersect at some point in the interior of the disk.  At this intersection point, the real and imaginary parts of $f(z)$ are both 0, so $f(z) = 0$; in other words, the intersection point is a root of $f$.

Gauss seemed to realize that his claim about algebraic curves had not been fully justified.  In a footnote, he wrote ``As far as I know, nobody has raised any doubts about this. However, should someone demand it then I will undertake to give a proof that is not subject to any doubt, on some other occasion.''  But in fact Gauss never gave such a proof.  In his footnote, Gauss went on to sketch a method of establishing the claim for the particular algebraic curves under consideration in his proof, but he didn't work out the details of this sketch.  The first exposition of Gauss's proof that included a complete justification for
Gauss's claim that any branch of the curve $\text{Re}(f(z)) = 0$ or $\text{Im}(f(z)) = 0$ that enters the disk $\abs{z} \le r$ must leave again
was given in 1920 by Alexander Ostrowski \cite{ostrowski}, and
this justification was not easy.
%Stillwell (p. 290) says:  "However, the existence of the connecting pieces is extremely hard to prove (and proving that they meet is not trivial either, being at least as hard as the intermediate value theorem).  The first proof was given by Ostrowski (1920).
(More recent versions of the proof can be found in \cite{gersten, martin, sjogren}.)  In discussing this point, Steve Smale wrote \cite[p.\ 4]{smale}:
%Smale also says (p. 5):  "I am grateful to Horst Simon for providing me with an English translation of Ostrowski's paper."

\begin{quote}
I wish to point out what an immense gap Gauss' proof contained.  It is a subtle point even today that a real algebraic curve cannot enter a disk without leaving.
\end{quote}

In the next section, we show that a small change in the strategy of the proof makes it possible to prove the existence of the required intersection point of the curves $\text{Re}(f(z)) = 0$ and $\text{Im}(f(z)) = 0$ without using anything more than elementary analysis.

\section{An elementary version of Gauss's proof.}

Following Gauss, we will prove the fundamental theorem for polynomials with real coefficients.  Suppose that $f$ is a polynomial of degree $N > 0$ with real coefficients.  By dividing by the leading coefficient, we may assume without loss of generality that $f$ is monic, so
\[
f(z) = z^N + \sum_{n=0}^{N-1} c_n z^n,
\]
where $c_0, \ldots, c_{N-1} \in \mathbb{R}$.  If $f(0) = 0$ then of course there is nothing to prove, so we may also assume $f(0) \ne 0$.

For real $r$ and $\theta$, we define
\begin{align*}
R_r(\theta) &= \text{Re}(f(re^{i\theta}))\\
&= \text{Re}\left[r^N(\cos\theta + i\sin\theta)^N + \sum_{n=0}^{N-1} c_n r^n(\cos\theta + i\sin\theta)^n\right],\\
I_r(\theta) &= \text{Im}(f(re^{i\theta}))\\
&= \text{Im}\left[r^N(\cos\theta + i\sin\theta)^N + \sum_{n=0}^{N-1} c_n r^n(\cos\theta + i\sin\theta)^n\right].
\end{align*}
Note that for fixed $r>0$, we can express each of these functions as a polynomial of degree $N$ in $\cos\theta$ and $\sin\theta$ by expanding the powers of $\cos\theta + i\sin\theta$.  Furthermore, the power of $\sin\theta$ in each term of $R_r(\theta)$ will be even, so by replacing $\sin^2\theta$ with $1-\cos^2\theta$ we can write $R_r(\theta)$ in the form
\[
R_r(\theta) = p_r(\cos\theta),
\]
where $p_r$ is a polynomial with real coefficients of degree at most $N$.  Similarly, the power of $\sin\theta$ in each term of $I_r(\theta)$ will be odd, so we have
\[
I_r(\theta) = \sin\theta\,q_r(\cos\theta),
\]
where $q_r$ is a polynomial with real coefficients of degree at most $N-1$.

From the formula $R_r(\theta) = p_r(\cos\theta)$ we see that the equation $R_r(\theta) = 0$ has at most $2N$ solutions (mod $2\pi$), and the only way it can have $2N$ solutions is if $p_r$ has $N$ distinct roots, all of which are in the interval $(-1,1)$, so that each is equal to $\cos\theta$ for two values of $\theta$ (mod $2\pi$).  In that case, $p_r$ factors into linear factors and its sign changes at each root, and it follows that the sign of $R_r$ changes at each zero.  In other words, if the zeros of $R_r$, listed in increasing order, are $\ldots, \theta_{-2}, \theta_{-1}, \theta_0, \theta_1, \theta_2, \ldots$, then the sign of $R_r(\theta)$ for $\theta_{j-1} < \theta < \theta_j$ is the opposite of the sign for $\theta_j < \theta < \theta_{j+1}$.

Similar conclusions hold for $I_r$.  From the formula $I_r(\theta) = \sin\theta\,q_r(\cos\theta)$ we see that every integer multiple of $\pi$ is a solution to the equation $I_r(\theta) = 0$.  The other solutions are values of $\theta$ for which $\cos\theta$ is a root of $q_r$.  There are at most $2N-2$ of these (mod $2\pi$), so there are at most $2N$ zeros of $I_r$ (mod $2\pi$), and if there are $2N$ zeros, then the sign of $I_r$ changes at each solution.

Next, we show that $R_r$ and $I_r$ have the maximum possible number of zeros when $r$ is large enough.  Intuitively, this follows from the fact that if $\abs{z}$ is large, then the $z^N$ term of $f(z)$ dominates the other terms, and the real and imaginary parts of $z^N$ change sign $2N$ times on any circle $\abs{z} = r$.  To make this idea precise, let
\[
r^* = \max\left(1,\sqrt{2}\sum_{n=0}^{N-1} \abs{c_n}\right),
\]
and consider any $r > r^*$.  For each integer $k$ in the range $0 \le k \le 4N$, let
\[
\theta_k = \frac{(2k-1)\pi}{4N}, \qquad z_k = re^{i\theta_k}.
\]
Then for every $k$,
\[
\abss{\sum_{n=0}^{N-1} c_n z_k^n} \le \sum_{n=0}^{N-1} \abs{c_n} r^n \le \left(\sum_{n=0}^{N-1} \abs{c_n}\right) \cdot r^{N-1} \le \frac{\sqrt{2}r^*}{2} \cdot r^{N-1} < \frac{\sqrt{2}}{2} \cdot r^N.
\]
Since $\arg(z_1^N) = N\theta_1 = \pi/4$,
\[
\text{Re}(z_1^N) = \frac{\sqrt{2}}{2} \cdot r^N > \abss{\sum_{n=0}^{N-1}c_nz_1^n}.
\]
Therefore $R_r( \theta_1) = \text{Re}(f(z_1)) > 0$.  On the other hand, $\arg(z_2^N) = 3\pi/4$, so
\[
\text{Re}(z_2^N) = -\frac{\sqrt{2}}{2} \cdot r^N < -\abss{\sum_{n=0}^{N-1} c_n z_2^n},
\]
and therefore $R_r( \theta_2) < 0$.  It follows that $R_r$ must have a zero in the interval $(\theta_1, \theta_2)$.  Similar reasoning can be used to show that $I_r( \theta_0) < 0 < I_r( \theta_1)$, so $I_r$ has a zero in the interval $(\theta_0, \theta_1)$.  Continuing to examine the signs of $R_r( \theta_k)$ and $I_r( \theta_k)$, we find that $R_r$ has a zero in each of the intervals $(\theta_1, \theta_2)$, $(\theta_3, \theta_4)$, \ldots, $(\theta_{4N-1}, \theta_{4N})$, and $I_r$ has a zero in each of the intervals $(\theta_0, \theta_1)$, $(\theta_2, \theta_3)$, \ldots, $(\theta_{4N-2}, \theta_{4N-1})$; the zero in the interval $(\theta_0, \theta_1) = (-\pi/(4N), \pi/(4N))$ is 0, and the zero in the interval $(\theta_{2N}, \theta_{2N+1}) = (\pi - \pi/(4N), \pi+ \pi/(4N))$ is $\pi$.  Thus the functions $R_r$ and $I_r$ both have $2N$ zeros (mod $2\pi$), and these zeros are interleaved: between any two zeros of either function is a zero of the other.

Let us say that a number $r > 0$ is an \emph{interleaving radius} if each of the functions $R_r$ and $I_r$ has $2N$ zeros (mod $2\pi$), and the zeros are interleaved.  Then what we have just shown is that all radii $r > r^*$ are interleaving.  For any interleaving $r$, let the zeros of $R_r$ in the interval $[0, 2\pi)$ be $\alpha_1(r), \ldots, \alpha_{2N}(r)$, and let the zeros of $I_r$ be $\beta_1(r), \ldots, \beta_{2N}(r)$, both listed in increasing order.  The fact that the roots are interleaved means that
\[
0 = \beta_1(r) < \alpha_1(r) < \beta_2(r) < \cdots < \beta_{2N}(r) < \alpha_{2N}(r) < 2\pi = \beta_1(r) + 2\pi.
\]
By our earlier observations, the signs of $R_r$ and $I_r$ change at each zero.

We claim now that the set of interleaving radii is open, and the functions $\alpha_j$ and $\beta_j$, for $1 \le j \le 2N$, are continuous on this set.  To see why this is true, suppose $c$ is interleaving and $\epsilon > 0$.  Choose a positive number $t \le \epsilon$ small enough that the intervals $(\alpha_j(c)-t, \alpha_j(c)+t)$ and $(\beta_j(c)-t,\beta_j(c)+t)$ are all disjoint (mod $2\pi$); in other words,
\begin{align*}
t = \beta_1(c)+t &< \alpha_1(c)-t,\\
\alpha_1(c)+t &< \beta_2(c)-t,\\
&\ \vdots\\
 \alpha_{2N}(c)+t &< 2\pi-t = \beta_1(c)-t+2\pi.
 \end{align*}
 Since the sign of $R_c$ changes at each zero, $R_c(\alpha_j(c)-t)$ and $R_c(\alpha_j(c)+t)$ must have opposite signs, and similarly $I_c(\beta_j(c)-t)$ and $I_c(\beta_j(c)+t)$ have opposite signs.  Now choose $\delta$ such that  $0 < \delta < c$ and $\delta$ is small enough that if $\abs{r-c} < \delta$, then for all $j$, $R_r(\alpha_j(c)-t)$ has the same sign as $R_c(\alpha_j(c)-t)$, and similarly for $R_r(\alpha_j(c)+t)$, $I_r(\beta_j(c)-t)$, and $I_r(\beta_j(c)+t)$.  Then $R_r(\alpha_j(c)-t)$ and $R_r(\alpha_j(c)+t)$ have opposite signs, and therefore $R_r$ has a zero in each of the intervals $(\alpha_j(c)-t, \alpha_j(c)+t)$.  Similarly, $I_r$ has a zero in each interval $(\beta_j(c)-t, \beta_j(c)+t)$.  Since these intervals are disjoint, it follows that $r$ is interleaving.  Note that the zero of $I_r$ in the interval $(\beta_1(c) - t, \beta_1(c) + t) = (-t, t)$ is $0 = \beta_1(r) = \beta_1(c)$.  Therefore $\alpha_j(r)$ belongs to the interval $(\alpha_j(c)-t, \alpha_j(c)+t)$, so $\abs{\alpha_j(r)-\alpha_j(c)} < t \le \epsilon$, and similarly $\abs{\beta_j(r)-\beta_j(c)} < \epsilon$.  This establishes that the set of interleaving radii is open, and the functions $\alpha_j$ and $\beta_j$ are continuous on this set.

\begin{figure}
\centering
\includegraphics[scale=.8]{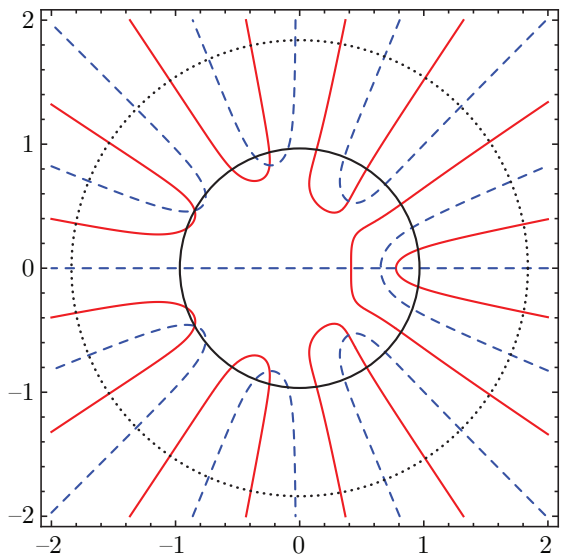}
\caption{The solid red lines are the points where $\text{Re}(f(z)) = 0$, and the dashed blue lines are the points where $\text{Im}(f(z)) = 0$, for the polynomial $f(z) = z^8+0.2z^7-0.1z^6-0.3z^5-0.1z^3+0.2z^2-0.3z+0.1$.  The large, dotted circle is $\abs{z} = r^*$, and the smaller, solid circle is $\abs{z} = r_0$.}\label{graph}
\end{figure}

Up to this point, our proof is not very different from Gauss's.  Our polar curves $\theta = \alpha_j(r)$ and $\theta = \beta_j(r)$ are parts of branches of the algebraic curves
$\text{Re}(f(z)) = 0$ and $\text{Im}(f(z)) = 0$
%$R_r( \theta) = 0$ and $I_r( \theta) = 0$ 
studied by Gauss.  Gauss's next step was to choose an interleaving radius $r$ and follow one of these branches from a point on the circle $\abs{z} = r$ into the interior of the circle, and he claimed that the branch would emerge at another point on the circle.  In contrast, our procedure is to follow all of the branches simultaneously from the circle $\abs{z} = r$ inward and show that two of the branches must eventually intersect.  This procedure is illustrated in Figure \ref{graph}.

Since $f(0) \ne 0$ and all coefficients of $f$ are real, $\text{Re}(f(0))$ is nonzero.  Therefore, by continuity, for all sufficiently small $r > 0$, $R_r$ has no zeros, so $r$ is not interleaving.  Thus the set of positive radii that are not interleaving is nonempty and bounded above, so it has a least upper bound $r_0$.  Since the set of interleaving radii is open, $r_0$ cannot be interleaving.

Let $(r_k)_{k=1}^\infty$ be a decreasing sequence of interleaving radii converging to $r_0$.  Applying the Bolzano--Weierstrass theorem to replace this sequence with a subsequence if necessary, we may assume that $\lim_{k \to \infty} \alpha_1(r_k)$ exists.  Passing to a further subsequence if necessary, we may also assume that $\lim_{k \to \infty} \beta_1(r_k)$ exists, and repeating this reasoning we may assume that all of the sequences $(\alpha_j(r_k))_{k=1}^\infty$ and $(\beta_j(r_k))_{k=1}^\infty$ converge.

Let $a_j = \lim_{k \to \infty} \alpha_j(r_k)$ and $b_j = \lim_{k \to \infty} \beta_j(r_k)$.  (It can in fact be shown that $a_j = \lim_{r \to r_0^+} \alpha_j(r)$ and $b_j = \lim_{r \to r_0^+} \beta_j(r)$, but we will not need these facts.)  By continuity of $R_r(\theta)$ and $I_r(\theta)$ (as functions of both $r$ and $\theta$), we have $R_{r_0}(a_j) = 0$ and $I_{r_0}(b_j) = 0$, and by the ordering of the functions $\alpha_j$ and $\beta_j$,
\[
0 = b_1 \le a_1 \le b_2 \le \cdots \le b_{2N} \le a_{2N} \le 2\pi = b_1 + 2\pi.
\]

If all of these inequalities are strict, then $r_0$ is interleaving, which is a contradiction.  Therefore at least one inequality is not strict.  In other words, for some $j$ and $k$, $a_j \equiv b_k \pmod{2\pi}$.  But then $R_{r_0}(a_j) = 0$ and $I_{r_0}(a_j) = I_{r_0}(b_k) = 0$, so $r_0e^{ia_j}$ is a root of $f$.

\section{Another approach.}

In the previous section, we have tried to stay as close as possible to Gauss's reasoning, in order to show how our proof gives a way to fill in the gap in Gauss's proof.  However, it turns out that the proof can be simplified a bit by deviating somewhat from Gauss's approach.  The main change is that, instead of intersecting the curves $\text{Re}(f(z)) = 0$ and $\text{Im}(f(z)) = 0$ with circles, we intersect them with horizontal lines.  In this version of the proof, there is no advantage to restricting attention to polynomials with real coefficients, so we prove the theorem directly for polynomials with complex coefficients.  Also, there is no need to introduce the cosine and sine functions into the proof; intersecting the curves $\text{Re}(f(z)) = 0$ and $\text{Im}(f(z)) = 0$ with lines rather than circles leads directly to polynomial equations.  We sketch this approach in this section, skipping those details that are essentially the same as in our first proof.

Let $f$ be a polynomial of degree $N > 0$ with complex coefficients.  Since the theorem is clearly true for linear polynomials, we may assume that $N \ge 2$, and as before we may also assume without loss of generality that $f$ is monic.  Thus, we can write
\[
f(z) = z^N + \sum_{n=0}^{N-1}c_n z^n,
\]
where $c_0, \ldots, c_{N-1} \in \mathbb{C}$.

For real numbers $x$ and $y$, we now define
\begin{align*}
R_y(x) &= \text{Re}(f(x+iy)),\\
I_y(x) &= \text{Im}(f(x+iy)).
\end{align*}
Straightforward algebra shows that we can write $R_y(x)$ and $I_y(x)$ in the form
\begin{align*}
R_y(x) &=  x^N + \sum_{n=0}^{N-1} g_n(y) x^n,\\
I_y(x) &= \sum_{n=0}^{N-1} h_n(y) x^n,
\end{align*}
where $g_n$ and $h_n$ are polynomials with real coefficients.  Furthermore, we have
\[
h_{N-1}(y) = Ny + a
\]
for some $a \in \mathbb{R}$.  Thus, for fixed $y$, $R_y$ is a monic polynomial of degree $N$ and $I_y$ is a polynomial of degree $N-1$, except that there is one value of $y$, namely $y = -a/N$, for which $I_y$ is either the constant zero polynomial or a nonzero polynomial of degree less than $N-1$.

We define a number $y$ to be \emph{interleaving} if the polynomial $R_y$ has exactly $N$ roots, $I_y$ has exacty $N-1$ roots, and the roots are interleaved.  In other words, if we let $\alpha_1(y), \ldots, \alpha_N(y)$ be the roots of $R_y$ and $\beta_1(y), \ldots, \beta_{N-1}(y)$ the roots of $I_y$, both listed in increasing order, then
\[
\alpha_1(y) < \beta_1(y) < \alpha_2(y) < \cdots < \beta_{N-1}(y) < \alpha_N(y).
\]

As in our previous proof, we can use the fact that the leading term of $f$ dominates the other terms to show that all sufficiently large $y$ are interleaving.  And as before, we can show that the set of interleaving values of $y$ is open, and the functions $\alpha_j$ and $\beta_j$ are continuous on this set.  However, not all values of $y$ are interleaving.  As we observed earlier, there is one value of $y$ for which $I_y$ is either the constant zero polynomial or a nonzero polynomial of degree less than $N-1$, and this value of $y$ cannot be interleaving.  Therefore the set of numbers that are not interleaving is nonempty, closed, and bounded above, so we can let $y_0$ be the largest number that is not interleaving.

The coefficients of $R_y$ are continuous functions of $y$, and are therefore bounded on $[y_0, y_0+1]$.  Using this fact, one can show that the function $\alpha_j$ and $\beta_j$ are bounded on $(y_0, y_0+1]$.  We can therefore use the Bolzano--Weierstrass theorem to find a decreasing sequence $(y_k)_{k=1}^\infty$ of interleaving numbers converging to $y_0$ such that the limits $a_j = \lim_{k \to \infty} \alpha_j(y_k)$ and $b_j = \lim_{k \to \infty} \beta_j(y_k)$ exist.  As before we have $R_{y_0}(a_j) = 0$, $I_{y_0}(b_j) = 0$, and
\begin{equation}\label{aborder}
a_1 \le b_1 \le a_2 \le \cdots \le b_{N-1} \le a_N.
\end{equation}

If $I_{y_0}$ is the constant zero polynomial, then for every $j$, $R_{y_0}(a_j) = 0$ and $I_{y_0}(a_j) = 0$, so $a_j + iy_0$ is a root of $f$.  Now suppose that $I_{y_0}$ is not the constant zero polynomial, so it has degree at most $N-1$.  If all of the inequalities in \eqref{aborder} are strict, then $y_0$ is interleaving, which is a contradiction.  Therefore at least one inequality is not strict.  In other words, for some $j$ and $k$, $a_j = b_k$.  But then $R_{y_0}(a_j) = 0$ and $I_{y_0}(a_j) = I_{y_0}(b_k) = 0$, so $a_j+iy_0$ is a root of $f$.

\bigskip
%\begin{acknowledgment}{Acknowledgment.}
\noindent\textbf{ACKNOWLEDGMENT.}
We would like to thank the anonymous referees for several helpful suggestions.
%\end{acknowledgment}

%\begin{biog}
%\item[Soham Basu] 
\noindent\textbf{Soham Basu}
received his B.Tech.\ from IIT Bombay in 2011 and dual M.Sc.\ from Imperial College London and TU Delft in 2013. Currently he is pursuing doctoral studies in photonics at \'Ecole Polytechnique F\'ed\'erale de Lausanne. He spends his spare time learning mathematics, music, and dance and photographing birds.
%\begin{affil}
%sohambasu6817@gmail.com
%\end{affil}

%\item[Daniel J. Velleman] 
\smallskip
\noindent\textbf{Daniel J. Velleman}
received his B.A.\ from Dartmouth
College in 1976 and his Ph.D.\ from the University of
Wisconsin--Madison in 1980.  He taught at the University of Texas
before joining the faculty of Amherst College in 1983.  Since 2011, he has also been an adjunct professor at the University of Vermont.  He was the editor of the American Mathematical Monthly from 2007 to 2011.  In his spare time he enjoys singing, bicycling, and playing volleyball.

%\begin{affil}
\noindent\textit{Department of Mathematics and Statistics,
Amherst College, Amherst, MA 01002\\
Department of Mathematics and Statistics, University of Vermont, Burlington, VT 05405}
%djvelleman@amherst.edu
%\end{affil}
%\end{biog}

\end{document}